\documentclass{amsart}

\usepackage{amssymb}
\usepackage{amsmath}
\usepackage{amsthm}

\def\E{\mathop{\hbox{\sf E}}\nolimits}
\def\P{\mathop{\hbox{\sf P}}\nolimits}

\def\Z{{\mathbb Z}}

\def\phi{\varphi}

\def\f{\frac}

\def\B{\Big}

\def\bs{\backslash}

\newtheorem{theorem}{Theorem}

\newtheorem{conjecture}{Conjecture}

\usepackage{graphics}

\def\eqref#1{(\ref{eq.#1})}

\def\defined#1{{\em #1}}

\def\optional#1{}

\def\proof{\par\noindent{\bf Proof.\ }}

\def\eop{\vskip 3mm }

\def\calR{{\cal R}}

\def\union{{\cup}}

\newtheorem{Lemma}{Lemma}

\def\cal{\mathcal}
\def\calQ{{\cal Q}}

\def\Sch{{\rm Sch}}
\def\Sk{{\rm Sk}}
\def\Sl{{\rm Sl}}
\def\Rg{{\rm Rg}}

\def\union{\cup}

\def\eop{\qed \vskip 3mm}

\begin{document}

\title[Distances in the infinite random quadrangulation]
      {On one property of distances in the infinite random quadrangulation}
\author{Maxim Krikun}
\address{Institut Elie Cartan Nancy, Nancy-Universit\'e, CNRS, INRIA, 
         Boulevard des Aiguillettes B.P. 239 F-54506 Vand{\oe}uvre l\`es Nancy}
\email{krikun@iecn.u-nancy.fr}
\thanks{Research supported by a grant BQR from R\'egion Lorraine.}

\begin{abstract}
We show that the Schaeffer's tree for an infinite quadrangulation
only changes locally when changing the root of the quadrangulation.
This follows from one property of distances in the infinite uniform
random quadrangulation. 
\end{abstract}

\maketitle

\section{Prerequisites}

\subsection{Quadrangulations and Schaeffer's bijection}

Let $Q$ be a rooted quadrangulation. Denote by 
$V(Q)$ the set of vertices of $Q$ and by $x$ the root vertex.
For every $z \in V(Q)$ let $d_x(z)$ be the (graph) distance from $z$ to $x$.

Quadrangulations are necessarily bipartite, 
so the increments of $d_x$ along any path in $Q$ are $\pm 1$.
Two vertices on the same distance from the root can not be joined by an edge,
but only by a diagonal in some face. Denote by $Q_\times$ the graph obtained by
adding to $Q$ the two diagonals of each face.
Note that $Q_\times$ is not a planar map anymore but only a graph.

The tree $\Sch_x(Q)$ is defined as follows \cite{ST}: 
for each face of $Q$ consider
the values of $d_x$ on the vertices of this face in counterclockwise order.
If for some $R$ these distances match the pattern $(R, R+1, R, R-1)$, 
add the edge $(R,R+1)$ to $\Sch_x$;
if the distances match the pattern $(R,R+1,R,R+1)$, add the diagonal $(R+1,R+1)$.
The graph thus obtained is a tree, spanning all the vertices of $Q$ except $x$.
In general $\Sch_x$ lives on the graph $Q_\times$, but not on the map $Q$ itself,
however, since $\Sch_x$ (as an embedded graph) has no self-intersections,
it inherits naturally a planar map structure from $Q$.

The main theorem concerning the tree $\Sch_x$ asserts that if for a 
finite sphere quadrangulation $Q$ one takes $\Sch_x$ (viewed as a planar tree)
and the values of $d_x$ on the vertices of $\Sch_x$, 
this information is sufficient to reconstruct $Q$. 
We omit the description of the actual reconstruction procedure, 
although it's very simple.

This bijection proved to be an important tool in the study of random quadrangulations,
and was used by many authors \cite{CD, MM, MM2, LG1, LG2, GM}.

On the other hand, the construction of $\Sch_x$ depends in an essential way on
the choice of the root in the quadrangulation. 
The main goal of this note is to gain a deeper understanding of this dependence.
In short, we show that this dependence is essentially continuous --- 
small displacement of the root in a typical random quadrangulation 
only leads to small perturbation in the structure of the tree.
To make this statement precise, we need to place ourselves in the framework of infinite
random maps.

\subsection{Infinite quadrangulation}
Consider a uniform measure $\mu_N$ on the set of rooted quadrangulations of a sphere
with $N$ faces. The local weak limit of this sequence
is a probability measure $\mu_\infty$ supported on infinite planar quadrangulations.
A sample from this measure is the \defined{uniform infinite planar quadrangulation}.

More precisely, let $\calQ$ be the space of all finite rooted quadrangulations, 
endowed with the ultrametric distance 
\[ D(q,q') = \inf\{ 2^{-R} | R\in\Z, B_R(q)=B_R(q') \}, \]
$B_R(q)$ denoting the ball or radius $R$ around the root in the quadrangulation $q$.
In this metric the sequence of quadrangulations $\{q_n\}_{n\in\Z}$ converges, 
if for every $R$ all of it's elements, starting from some index $n(R)$,
coincide on a ball of radius $R$ around the root.
Elements of the completion $\bar\calQ$, different from finite quadrangulations,
are called \defined{infinite quadrangulations}, and the measure $\mu_\infty$ is 
supported on $\bar\calQ \bs \calQ$.
This construction was first considered in~\cite{AS} in the case of triangulations;
the treatment for quadrangulations in completely analogous~\cite{K}.

\subsection{Invariance with respect to the random walk}\label{sec.inv}
Let $Q$ be a finite quadrangulation, and let $\{v_t\}_{t\in\Z}$ be the simple 
random walk on $Q$ (viewed as a graph). Let $e_t$ be the edge traversed by a random
walk prior to arriving at $v_t$, then $\{e_t\}_{t\in\Z}$ is a random walk on the
set of directed edges of $Q$ (the \defined{edge random walk}),
and it's stationary measure is uniform on this set.

Therefore, if one acts upon the space of probability measures on finite quadrangulations
by displacing the root randomly according to one-step transitions of the edge random walk,
the measures $\mu_N$ are invariant under this action, 
and so is the limiting measure $\mu_\infty$.

From this observation we conclude, in particular, that if some property holds with
probability one for a uniform infinite quadrangulation, this property will hold 
with probability one for the same quadrangulation with the root displaced for
a uniformly bounded distance according to any rule, deterministic or random.

\subsection{Topological ends}
An infinite connected graph $G$ has \defined{one end}
if for any finite subgraph $H\subset G$ at most one (and therefore exactly one) 
component of the complement $G\bs H$ is infinite.
In particular, if $G$ is a tree then from each vertex of $G$
there is exactly one infinite simple path. 
If $G$ is planar and has one end, then it's planar dual also has one end.

Both the uniform infinite triangulation and the uniform infinite 
quadrangulation have one end \cite{AS, K}.

\subsection{Skeleton of the infinite quadrangulation}\label{sec.skel}
Given an infinite rooted quadrangulation $Q$, consider the set of vertices 
situated at distance $R$ from the root and the edges of $Q_\times$
between these vertices. Denote the graph thus obtained by $\Gamma_R$.
$\Gamma_R$ cuts $Q$ into a number of connected components, only one 
of which is infinite a.s.~(this follows from the a.s.~one-endedness of $Q$).
Let $\gamma_R$ be a cycle in $\Gamma_R$ that separates the root from
this infinite component;
if some vertices of $\Gamma_R$ are connected by more than one edge,
and such a cycle happens to be non-unique, 
let us choose among all alternatives the edge which is closer to the root.

Each edge of $\gamma_R$ is an $(R,R)$-diagonal in some face of type $(R, R-1, R, R+1)$, 
and the vertex corresponding to $(R-1)$ belongs necessarily to the cycle $\gamma_{R-1}$.
Let $e_1, e_2, \ldots$ be the edges of $\gamma_R$, 
and $w_1, w_2, \ldots$ the associated vertices on $\gamma_{R-1}$.

We introduce a directed tree structure on the union of edges of $\gamma_R$, 
by declaring the edges between $w_n$ and $w_{n+1}$
to be the offspring of the edge $e_n$.
The tree thus obtained is called the \defined{skeleton} of the quadrangulation $Q$
and is written $\Sk(Q)$ or $\Sk_x(Q)$ to emphasise the dependence on the choice of the root.
Note that $\Sk(Q)$ is directed from infinity to the root.

When $Q$ is the uniform infinite quadrangulation,
$\Sk(Q)$ is distributed as the tree of a certain time-reversed critical Galton-Watson process,
whose offspring distribution can be calculated explicitly \cite{K}.
%
In particular it follows that for every $R$ with probability one 
all of the individuals on $\gamma_R$ have a single common ancestor, 
and therefore tree $\Sk(Q)$ has one end (see \ref{sec.calc1} for details).

\section{Main theorem}

Let $x$ and $y$ be two distinct vertices of the quadrangulation $Q$;
we may assume that $x$ is the root vertex.
Let $f(z) = d_x(z) - d_y(z)$ be a function on the vertices of $V(Q)$.

From the definition of $\Sch_x$ it follows that for each face of $Q$ 
the structure of $\Sch_x$ inside it only depends 
on the increments of $d_x$ around this face.
Therefore if the quadrangulation $Q$ is infinite
and $f$ is constant in $Q$ except for a finite number of vertices,
then the trees $\Sch_x$ and $\Sch_y$ coincide almost everywhere in $Q$.

For simplicity let's take $y$ to be a neighbour of $x$, i.e.~$d_x(y)=1$.
Then we can formulate the main result as follows:
\begin{theorem}
Let $(x,y)$ be the root in a uniform infinite random quadrangulation $Q$.
With probability one 
the function $f(z) = d_x(z) - d_y(z)$ is constant almost everywhere in $Q$.
\end{theorem}

\section{Proof}

\subsection{Trees and geodesics within quadrangulations}

A \defined{geodesic emanating from $x$} is a finite or infinite path 
starting at $x$ and such that $d_x$ is strictly increasing along this path.

Let $z\in V(Q)$. An \defined{$x$-slit from $z$} is an infinite simple 
path in $Q_\times$ starting at $z$,
such that $d_x$ is non-decreasing along this path.
Denote by $\Sl_x(Q)$ the set of vertices of $Q$ admitting an $x$-slit.

Note that since the number of vertices at any finite distance from $x$ is 
finite a.s., $d_x$ tends to infinity along any $x$-slit.
Also note that the $x$-slit does not necessarily exist for every $z\in V(Q)$.

Given a point $z \in \Sl_x(Q)$, we define the \defined{rightmost geodesic} 
from $z$ to $x$, denoted $\Rg^z_x$, recursively as following:
let $a$ be the first point on the $x$-slit from $z$, and let $b$ be the first
point after $A$ in the clockwise order around $z$, such that $d_x(b) = d_x(z)-1$.
Then $\Rg^z_x$ consists of the edge $(z,b)$ followed by the rightmost geodesic 
from $b$ to $x$.

Clearly, the choice the vertex $b$ does not depend on the choice of the $x$-slit;
and given the $x$-slit from $z$ exists, an $x$-slit from $b$ exists as well,
so the rightmost geodesic from $b$ is well-defined. 

We call an infinite geodesic $\gamma$, emanating from $x$, \defined{rightmost}, 
if for every $z\in\gamma$ it is the rightmost geodesic from $z$ to $x$.

All rightmost geodesics on $Q$ form a tree; denote this tree by $\calR$.
This tree can be seen as a dual tree to $\Sk(Q)$, and therefore it has one end.
Thus we obtain 
\begin{Lemma}\label{L1}
The rightmost infinite geodesic $\Rg_X^\infty$ is unique a.s.
\end{Lemma}

Now let $y$ be some vertex in $Q$ different from $x$, and let $f(z) = d_x(z) - d_y(z)$.
For simplicity we'll take $d_x(y)=1$, then $f$ only takes two values $\pm1$.
There are two alternatives -- 
either $f$ takes the same value for all $z\in V(Q)$ except a finite subset, 
either each of the two sets $V^\pm := \{ z \in V(Q) | f(z) = \pm1\}$ is infinite. 
Our goal is to prove that the second case doesn't happen,
therefore from now on we assume that both $V^+$ and $V^-$ are infinite.

\begin{Lemma}\label{L2}
$f$ is non-decreasing along any geodesic emanating from $x$;
$f$ is non-increasing along any geodesic emanating from $y$.
\end{Lemma}

\proof Let $\gamma$ be a geodesic emanating from $x$. The increments of $d_x$
along $\gamma$ are all equal to $+1$, while the increments of $d_y$ are $\pm1$, 
so the statement follows. The statement for $y$ follows by symmetry.
\eop

\begin{Lemma}\label{L2a}
$f$ is constant almost everywhere in~$V(Q)$ if and only if $f$ is constant on $\gamma_R$ for some $R > d(x,y)$.
\end{Lemma}

\proof If $f$ is constant on $V(Q)$ except the finite set $W$, then $f$ is constant on 
$\gamma_R$ for all $R > {\rm diam}(W)$.

Suppose now that $f$ is constant on some cycle $\gamma_R$, $R > d(x,y)$.
As mentioned above, $\gamma_R$ separates $Q$ in two parts,
and the part $W$ containing the root $x$ is finite~a.s.
Since $R>d(x,y)$, $W$ also contains $y$.
Pick a point $z \notin W$ and some geodesic from $z$ to $x$; this geodesic intersects
$\gamma_R$ at some point $x'$. Similarly, take some geodesic from $z$ to $y$ and let $y'$ 
be the point where this geodesic intersects $\gamma_R$.

By the previous lemma, $f$ is non-decreasing on the geodesic from $x'$ to $z$,
and non-increasing on the geodesic from $y'$ to $z$, 
and since $f$ is constant on $\gamma_R$ we conclude that $f(x') = f(y') = f(z)$,
thus $f$ is constant on the complement of $W$.
\eop

\begin{Lemma}\label{L3}
Assume that both $V^-$ and $V^+$ are infinite. 
Then $f=-1$ on $\Rg_x^\infty$.
\end{Lemma}
\proof 
Assume the opposite, i.e.~that $f=+1$ on $\Rg_x^\infty$ starting from some point $v\in\Rg_x^\infty$.
Since $f$ is nondecreasing along each geodesic emanating from $x$,
$f=+1$ on the whole subtree of $\calR$ above $v$.
But since $\calR$ is locally finite and has one end a.s.,
it follows that $f=-1$ for at most a finite number of vertices of $\calR$.
Therefore $f=+1$ on $\gamma_R$ for some large $R$ and by the previous lemma
$f$ is constant a.e.~in $V(Q)$ --- a contradiction.
\eop

%
%
%

Because of the invariance with respect to the random walk,
every property that holds a.s.~for $Q$ rooted at $x$
holds also for $Q$ rooted at $y$.
Using this, and the fact that $f$ changes it's sign if we swap $x$ with $y$,
we obtain the following
\begin{Lemma}\label{L4}
There exists an infinite rightmost geodesic $\Rg_y^\infty$,
and $f=+1$ on this geodesic unless $V^+$ is finite.
\end{Lemma}
Note that the rightmost geodesic to $y$ should be considered with respect
to the $y$-slits, and as an $x$-slit is not necessarily a $y$-slit, 
$\Sl_x(Q)$ does not necessarily coincide with $\Sl_y(Q)$.

\subsection{}
Now assuming both $V^+$ and $V^-$ are infinite,
we have the following picture: there exists a rightmost geodesic for $x$, 
and a rightmost geodesic for $y$, and these two are disjoint.

Observe that $\Rg_x^\infty$ is a geodesic for $y$: since $f$ is constant,
and $d_x$ is strictly increasing along $\Rg_x^\infty$, $d_y$ is strictly 
increasing as well, so adding the edge $(x,y)$ to $\Rg_x^\infty$ we obtain
an honest infinite geodesic emanating from $y$.

By Lemma~\ref{L4} the rightmost infinite geodesic for $y$ is unique,
so $\Rg_x^\infty$ should not be rightmost for $y$. 
More precisely, for every point $z\in \Rg_x^\infty$ 
define the \defined{divergence point} $\delta(z)$ as the lowest
common point (lowest in the sense of $d_x$ and $d_y$) of the geodesics $\Rg_y^z$ and $\Rg_x^z$.
Then the set of divergence points
\[ \Delta := \B\{ v \,\B|\, \mbox{$v=\delta(z)$ for some $z\in \Rg_x^\infty$} \B\} \] 
has to be infinite.

Let $v\in\Delta$ be a divergence point, and let $R=d_x(v)$.
Let $b$ be the first point on the geodesic $\Rg_y^v$.
Then we have 
\begin{itemize}
\item $d_y(b) = d_y(v)-1 = R$, 
\item $d_x(b) > d_x(v)-1$ therefore $d_x(b) = R+1$,
\item $f(b)=+1$.
\end{itemize}

Let's now compare the rightmost geodesic $\Rg_y^b$ with the rightmost geodesic $\Rg_x^b$.
$\Rg_y^b \union (y,x)$ is a geodesic from $b$ to $x$, therefore $\Rg_x^b$ 
either coincides with $\Rg_y^b$ or diverges from it and goes to the right\footnote{
i.e.~to the left, from the point of view of an observer moving along the geodesic
towards $y$.}

Note that $\Rg_x^b$ and $\Rg_y^b$ may coincide only for a finite number 
of divergence points $z$ --- otherwise it would follows that $\Rg_y^\infty$ 
is a second infinite rightmost geodesic for $x$, which as we know does not happen.

On the other hand, it's clear that $\Rg_x^b$ and $\Rg_y^b$
may not diverge until $\Rg_y^b$ merges with $\Rg_y^\infty$
(otherwise $\Rg_x^b$ would cross $\Rg_y^\infty$ in some point, and $\Rg_y^b$ 
would not be the rightmost geodesic from $b$ to $y$).

Therefore (for a.e.~divergence point in $\Delta$) the geodesic $\Rg_x^b$ 
crosses $\Rg_y^\infty$, and eventually merges with $\Rg_x^\infty$ 
but on the other side, so that parts of the geodesics~$\Rg_x^v$, $\Rg_x^b$ 
together with the edge $(v,b)$ form a loop that separates the root $x$ from infinity.
It follows that the trunk of the skeleton $\Sk_x(Q)$ passes between $v$ and $b$.

Since $v \in \gamma_R$, we obtain that the trunk of $\Sk_x(Q)$ intersects $\gamma_R$
at an edge immediately next to $v$, 
and $v$ is the point where $\Rg_x^\infty$ intersects $\gamma_R$.

To conclude the proof of the main theorem we will show that this can only 
happen finitely many times.

\begin{Lemma}\label{L5}
For a given $R$ let $\gamma_R$ be the cycle in $Q_\times$ that is situated at distance
$R$ from the root and that separates the root from the infinite part of $Q$.
Let $v_R$ be the vertex where $\Rg_x^\infty$ intersects $\gamma_R$,
and let $e_R$ be edge where the trunk of $\Sk_x(Q)$ intersects $\gamma_R$.

Then for any $d_0$ with probability one $d(v_R,e_R) < d_0$ only for finitely many $R$.
\end{Lemma}

\proof
The cycle $\gamma_R$ separates the quadrangulation into two parts,
and it's known \cite{K} that given the length of $\gamma_R$
these parts are conditionally independent.

Observe that the position of $e_R$ on $\gamma_R$ only depends on the 
finite part (the one containing $x$), while the position of $v_R$ only
depends on the infinite part. 
The measure $\mu_\infty$ is invariant with respect to cutting the quadrangulation 
along $\gamma_R$ and gluing the parts back with a uniformly distributed twist.
Due to such rotational symmetry, the probability that $e_R$ and $v_R$ are close
is of order $1/|\gamma_R|$.

Since $|\gamma_R|$ grows as $R^2$ (and it can be checked, see~\ref{sec.calc2}, that $\E |\gamma_R|^{-1} = O(R^{-2})$)
it follows by Borel-Cantelli that only a finite number of events
$\{d(v_R, e_R) < d_0\}_{R=0}^\infty$ can occur simultaneously.
\eop
This finishes the proof of the main theorem.


\section{Concluding remarks}

For every directed edge $(x,y)$ in $Q$ define
\[ \Delta h(x,y) = \lim_{z\to\infty} d(x,z) - d(y,z), \]
then there exists a function $h$ on $V(q)$, defined up to an additive constant by
%
\[ \Delta h(x,y) = h(x) - h(y). \]
Applying the construction of the tree $\Sch$ with repsect to the increments
of $h$, one obtains a spanning forest  of $Q_\times$. Denote this graph by $\Sch_\infty$.
It seems natural to interpret $\Sch_\infty$ as a limit of trees $\Sch_z(Q)$, when $z$ is sent 
to infinity, and to conjecture that
\begin{conjecture}
$\Sch_\infty(Q)$ is connected a.s.
\end{conjecture}

From a slightly different point of view, one can consider $\Sch_\infty$ as 
the (limit of the) Schaeffer's tree in a large random quadrangulation, 
seen from a uniformly chosen point.
\begin{conjecture}
Let $Q_n$ be a uniform rooted quadrangulation with $n$ faces, and $z_n$ a uniformly 
chosen vertex in $Q_n$. Then the pair $(Q_n, \Sch_{z_n}(Q_n))$ converges as $n\to\infty$
to $(Q, \Sch_\infty)$ (in the sense of local weak convergence).
\end{conjecture}

Provided that the first conjecture above holds, it would be useful to know the distribution
of the tree $\Sch_\infty(Q)$, seen as a planar tree with~a~${-1,0,1}$-valued antisymmetric 
function $\Delta h$ on directed edges. The representation of the tree $\Sch_x(Q)$ in \cite{CD}
strongly suggests the following
\begin{conjecture}
The pair $(\Sch_\infty, \Delta h)$ is distributed as a Galton-Watson tree with 
geometric-$1/2$ offspring distribution, conditioned to be infinite, and with increments
distributed independently and uniformly among $\{-1,0,1\}$.
\end{conjecture}




%

\newpage

\section{Some additional calculations}

\subsection{}\label{sec.calc1}
Let $Q$ be a uniform infinite quadrangulation,
and $\gamma_R$ be a sequence of cycles in $Q_\times$, 
located at distance $R$ from the root and separating the root 
from infinity, as defined in Section~\ref{sec.skel}.
In~\cite{K} the following theorem was proved:
\begin{theorem}\label{T2}
$|\gamma_R|$ is a Markov chain with transition probabilities given by
\[ \P\B\{ |\gamma_{r+n}|=k \B| |\gamma_{r}|=l \B\} 
    = \f{[t^k]F(t)}{[t^l]F(t)} \cdot \P\{ \xi_n = l | \xi_0 = k \},
\]
where $\xi$ is a critical branching process with offspring generating function
\[ \phi(t) = \f1{2t} \B( \sqrt{(t-9)(t-1)^3} - 3 + 6t - t^2\B), \]
and 
\[ F(t) = \f34 \B( \sqrt{\f{9-t}{1-t}} - 3 \B). \]
\end{theorem}

Moreover, the tree structure $\Sk(Q)$ restricted to $\gamma_1, \ldots, \gamma_R$,
and conditioned on $|\gamma_R| = m$ has the distribution of a 
forest of $m$ Galton-Watson trees, 
conditioned to have exactly one individual at height $R$
and no individuals higher than $R$.

Let $A_R^h$ be the number individuals at $\gamma(R+h)$, having nonempty offspring at $\gamma(R)$.
Then the generating function for $A_R^{R+h}$ is given by
\[ \E y^{A_R^{R+h}} = [t] F\B( \phi_h(0) + y (\phi_{R+h}(t) - \phi_h(0)) \B), \]
where $\phi_R$ is the $R$-fold iteration of $\phi$, 
which can be written explicitly as
\begin{equation}\label{eq.phiiter}
 \phi_R(t) = 1 - \f8{ \B(\sqrt{\f{9-t}{1-t}}+2R \B)^2 - 1}. 
\end{equation}
In particular we obtain 
\begin{eqnarray*}
\P\{ A_R^h = 1 \} &=&
[y] [t] F\B( \phi_h(0) + y (\phi_{R+h}(t) - \phi_h(0)) \B) \\
&=& \f{ (h+2)^2(h+1)^2 (2R + 2h +3) }{ (2h+3)(R+h+2)^2 (R+h+1)^2 }.
\end{eqnarray*}
It follows that
\[ \lim_{h\to\infty} \P \{ A_R^h = 1\} = 1, \]
so all individuals at $\gamma_R$ have a single common ancestor a.s,
thus proving that the graph $\Sk(Q)$ is connected a.s.,
and
\[ \sum_{h\ge0} \P \{ A_R^h = 1 \} = \infty, \]
so the distance to this common ancestor has infinite mean.

\subsection{}\label{sec.calc2}
Let
\[ \Phi(t) = \int_0^t \f{F(y)}{y} dy, \]
then
\[ \E |\gamma_R|^{-1} = [t] \Phi(\phi_R(t)) = \f{2(2R+3)}{(R+1)(R+2)(R+3)} = O(R^{-2}). \]

\end{document}